\documentclass[10pt]{elsarticle}
\journal{}
\usepackage{amsmath}
\usepackage{indentfirst}
\usepackage{amsfonts,amssymb}
\usepackage{mathrsfs}
\usepackage{amsthm}
\usepackage{lineno,hyperref}
\modulolinenumbers[5]
\newtheorem{Definition}{Definition}[section]
\newtheorem{Theorem}{Theorem}[section]
\newtheorem{Lemma}{Lemma}[section]

\numberwithin{equation}{section}

\begin{document}
	
	\begin{frontmatter}
		
		\title{Dynamics of an imprecise stochastic Holling II one-predator two-prey
			system with jumps}

		\author{Fei Sun\corref{mycorrespondingauthor}}
		\cortext[mycorrespondingauthor]{Corresponding author}
		\address{School of Mathematics and Computational Science, Wuyi University, Jiangmen 529020, China}
		\ead{sunfei@whu.edu.cn (fsun.sci@outlook.com)}

\begin{abstract} Groups in ecology are often affected by sudden environmental perturbations. Parameters of stochastic models are often imprecise due to various uncertainties. In this paper, we formulate a stochastic Holling II one-predator two-prey system with jumps and interval parameters. Firstly, we prove the existence and uniqueness of the positive solution. Moreover, the sufficient conditions for the extinction and persistence in the mean of the solution are obtained. %Finally, some simulations are carried out to demonstrate our theoretical results.
\end{abstract}


\end{frontmatter}

\section{Introduction}
\label{sec:1}

In ecology and mathematical ecology, the study of interrelationship between species has become one of the main topics. And there have been growing interests on the dynamical behavior of the population species living in groups, such as Holling type I, II, and III functional response. 
 For a better review of Holling II functional response and its extension, see \cite{1}- \cite{9} as well as there references. 

However, the sudden environmental perturbations may bring substantial social and economic losses. For example, the recent COVID-19 has a serious impact on the world. It is more realistic to study the population dynamics with imprecise parameters. Panja et al. \cite{23} studied a cholera epidemic model with imprecise numbers and discussed the stability condition of equilibrium points of the system. Das and Pal \cite{24} analyzed the stability of the system and solved the optimal control problem by introducing an imprecise SIR model. Other studies on imprecise parameters include those of \cite{20}-\cite{222}, and the references therein.

The main focus of this paper is dynamics of an imprecise stochastic Holling II one-predator two-prey
model with jumps. To this end, we first introduce the imprecise stochastic Holling II one-predator two-prey model. With the help of Lyapunov functions, we prove the existence and uniqueness of the positive solution. Further, the sufficient conditions for the extinction and persistence in the mean of the solution are obtained.

The remainder of this paper is organized as follows. In Sect.~\ref{sec:2}, we introduce the basic models. In Sect.~\ref{sec:3}, the unique global positive solution of the system is proved. 
The sufficient conditions for the extinction and persistence in the mean of the solution are derived in
Sect.~\ref{sec:4}. %Finally, in Sect.~\ref{sec:5}, numerical simulations and discussions are presented. \\

\section{Imprecise stochastic Holling II one-predator two-prey system}
\label{sec:2}

In this section, we introduce the imprecise stochastic system.
Let $ x_{i}(t) $ $ (i=1,2) $ and $ y(t) $ denote the population sizes of prey species  and the population size of predator species at time $ t $, respectively.  Then a stochastic Holling II one-predator two-prey system takes the following form \cite{27}.

\begin{equation}\label{2.1}
\left\{
\begin{array}{lcl}
dx_{1}(t) = x_{1}(t) [r_{1} - a_{11}x_{1}(t) - a_{12}x_{2}(t) - \dfrac{a_{13}y(t)}{1+x_{1}(t)}]dt + \sigma_{1}x_{1}(t)dB_{1}(t) +  \int_{\mathbb{Y}}c_{1}(u)x_{1}(t^{-}){N}(dt,du),\\
dx_{2}(t) = x_{2}(t) [r_{2} - a_{21}x_{1}(t) - a_{22}x_{2}(t) - \dfrac{a_{23}y(t)}{1+x_{2}(t)}]dt + \sigma_{2}x_{2}(t)dB_{2}(t) +  \int_{\mathbb{Y}}c_{2}(u)x_{2}(t^{-}){N}(dt,du),\\
dy(t) = y(t) [-r_{3} - a_{33}y(t) + \dfrac{a_{31}x_{1}(t)}{1+x_{1}(t)} + \dfrac{a_{32}x_{2}(t)}{1+x_{2}(t)}]dt + \sigma_{3}y(t)dB_{3}(t) +  \int_{\mathbb{Y}}c_{3}(u)y(t^{-}){N}(dt,du),
\end{array}  
\right.
\end{equation}
where $  x_{1}(t^{-}) $, $ x_{2}(t^{-}) $ and $ y(t^{-}) $ are the left limits of  $ x_{1}(t) $, $ x_{2}(t) $ and $ y(t) $, respectively. $ r_{i}>0 $ $ (i=1,2,3) $ are the intrinsic growth rates or death rate, $ a_{ii} $ $ (i=1,2,3) $ stand for the intraspecies interaction,  $ a_{ij} $ $ (i\neq j) $ represent the effect of species $ j $ upon the growth rate of species $ i $.
$ B_{i}(t) $, $( i = 1, 2, 3) $ are mutually independent Brownian motion defined on a complete probability space $ (\Omega, \mathcal{F}, \mathcal{F}_{t\geq 0}, \mathbb{P}) $.
 $ \sigma_{i}^{2} $ represent the intensities of  $ B_{i}(t) $.
  Let $ \lambda $ be the characteristic measure of $ N $ which is defined on a finite measurable subset $ \mathbb{Y} $ of $ (0,+\infty) $ with $ \lambda(\mathbb{Y}) < \infty $. 
Define the compensated random measure by $ \widetilde{N}(dt, du) = N(dt, du) − \lambda(du)dt $. \\

%We also denote $ \mathbb{R}^{d}_{+} = \{x \in \mathbb{R}^{d} : x_{i} > 0 \textrm{ for all } 1 \leq i \leq d\} $, $ \overline{\mathbb{R}}^{d}_{+} = \{x \in \mathbb{R}^{d}: x_{i} \geq 0 \textrm{ for all } 1 \leq i \leq d\} $. If $ f(t) $ is an integrable function on $  [0, \infty) $, define $\langle f \rangle_{t} =\dfrac{1}{t} \int_{0}^{t} f(s) ds$.\\

Before we state the imprecise stochastic Holling II one-predator two-prey system, definitions of Interval-valued function should recalled (Pal \cite{24}).

\begin{Definition}
(Interval number) An interval number $ A $ is represented by closed interval $ [g^{l}, g^{u}] $ and defined by $ A = [g^{l}, g^{u}]  = \{x|g^{l} \leq x \leq g^{u}, x \in \mathbb{R}\},$ where $ \mathbb{R} $ is the set of real numbers and $ g^{l} $, $  g^{u} $ are the lower and upper limits of the interval numbers, respectively. The interval number $ [g, g] $ represents a real number $ g $. The arithmetic operations for any two interval numbers $ A = [g^{l}, g^{u}]$ and $ B = [h^{l}, h^{u}]$ are as follows:\\
\indent Addition: $ A+B= [g^{l}, g^{u}] + [h^{l}, h^{u}] = [g^{l}+h^{l}, g^{u}+h^{u}] $.\\
\indent Subtraction:$ A-B= [g^{l}, g^{u}] - [h^{l}, h^{u}] = [g^{l}-h^{l}, g^{u}-h^{u}] $.\\
\indent Scalar multiplication: $\alpha A= \alpha [g^{l}, g^{u}] =  [\alpha g^{l}, \alpha g^{u}]$, where $ \alpha $ is a positive real number.\\
\indent  Multiplication: $ AB= [g^{l}, g^{u}]  [h^{l}, h^{u}]= [\min \{g^{l}h^{l}, g^{u}h^{l}, g^{l}h^{u}, g^{u}h^{u} \}, \max \{g^{l}h^{l}, g^{u}h^{l}, g^{l}h^{u}, g^{u}h^{u} \}  ]  $.\\
\indent  Division: $ A/B =[g^{l}, g^{u}]/[h^{l}, h^{u}]=[g^{l}, g^{u}][ \frac{1}{h^{l}}, \frac{1}{h^{u}}] $.
\end{Definition}

\begin{Definition}
(Interval-valued function) Let $ g > 0 $, $ h > 0 $. If the interval is of the from $ [g, h] $, the interval-valued function is take as $ f(k) = g^{(1-k)}h^{k} $ for $ k \in [0, 1] $.	
\end{Definition}

Let $ \hat{r}_{i}, \hat{a}_{ij}, \hat{\sigma}_{i}$ represent the interval numbers of $ r_{i}, {a_{ij}}, {\sigma_{i}}$ $ (i,j=1,2,3) $, respectively. The system (\ref{2.1}) with imprecise parameters becomes:

\begin{equation}\label{2.2}
\left\{
\begin{array}{lcl}
dx_{1}(t) = x_{1}(t) [\hat{r}_{1} - \hat{a}_{11}x_{1}(t) - \hat{a}_{12}x_{2}(t) - \dfrac{\hat{a}_{13}y(t)}{1+x_{1}(t)}]dt + \hat{\sigma}_{1}x_{1}(t)dB_{1}(t) +  \int_{\mathbb{Y}}c_{1}(u)x_{1}(t^{-}){N}(dt,du),\\
dx_{2}(t) = x_{2}(t) [\hat{r}_{2} - \hat{a}_{21}x_{1}(t) - \hat{a}_{22}x_{2}(t) - \dfrac{\hat{a}_{23}y(t)}{1+x_{2}(t)}]dt + \hat{\sigma}_{2}x_{2}(t)dB_{2}(t) +  \int_{\mathbb{Y}}c_{2}(u)x_{2}(t^{-}){N}(dt,du),\\
dy(t) = y(t) [-\hat{r}_{3} - \hat{a}_{33}y(t) + \dfrac{\hat{a}_{31}x_{1}(t)}{1+x_{1}(t)} + \dfrac{\hat{a}_{32}x_{2}(t)}{1+x_{2}(t)}]dt + \hat{\sigma}_{3}y(t)dB_{3}(t) +  \int_{\mathbb{Y}}c_{3}(u)y(t^{-}){N}(dt,du),
\end{array}  
\right.
\end{equation}

where $ \hat{r}_{i}= [r_{i}^{l}, r_{i}^{u}] $, $ \hat{a}_{ij}= [a_{ij}^{l}, a_{ij}^{u}] $, $ \hat{\sigma_{i}}= [\sigma_{i}^{l}, \sigma_{i}^{u}] $ $ (i,j=1,2,3) $.  \\

According to the Theorem 1 in Pal et al. \cite{20} and considering the interval-valued function $ f (p) = (f^{l})^{1-p}(f^{u})^{p} $ for interval $ \hat{f}= [f^{l}, f^{u}] $  for $ p \in [0, 1] $, we can prove that system (\ref{2.2}) is equivalent to the following system:

\begin{equation}\label{2.3}
\left\{
\begin{array}{lcl}
dx_{1}(t) = x_{1}(t) [(r_{1}^{l})^{1-p}(r_{1}^{u})^{p} - (a_{11}^{l})^{1-p}(a_{11}^{u})^{p}x_{1}(t) - (a_{12}^{l})^{1-p}(a_{12}^{u})^{p}x_{2}(t) - \dfrac{(a_{13}^{l})^{1-p}(a_{13}^{u})^{p}y(t)}{1+x_{1}(t)}]dt\\
 \ \ \ \ \ \ \ \ \ \ \ \  + (\sigma_{1}^{l})^{1-p}(\sigma_{1}^{u})^{p}x_{1}(t)dB_{1}(t) +  \int_{\mathbb{Y}}c_{1}(u)x_{1}(t^{-}){N}(dt,du),\\
dx_{2}(t) = x_{2}(t) [(r_{2}^{l})^{1-p}(r_{2}^{u})^{p} - (a_{21}^{l})^{1-p}(a_{21}^{u})^{p}x_{1}(t) - (a_{22}^{l})^{1-p}(a_{22}^{u})^{p}x_{2}(t) - \dfrac{(a_{23}^{l})^{1-p}(a_{23}^{u})^{p}y(t)}{1+x_{2}(t)}]dt\\
 \ \ \ \ \ \ \ \ \ \ \ \ + (\sigma_{2}^{l})^{1-p}(\sigma_{2}^{u})^{p}x_{2}(t)dB_{2}(t) +  \int_{\mathbb{Y}}c_{2}(u)x_{2}(t^{-}){N}(dt,du),\\
dy(t) = y(t) [-(r_{3}^{l})^{1-p}(r_{3}^{u})^{p} - (a_{33}^{l})^{1-p}(a_{33}^{u})^{p}y(t) + \dfrac{(a_{31}^{l})^{1-p}(a_{31}^{u})^{p}x_{1}(t)}{1+x_{1}(t)} + \dfrac{(a_{32}^{l})^{1-p}(a_{32}^{u})^{p}x_{2}(t)}{1+x_{2}(t)}]dt\\
  \ \ \ \ \ \ \ \ \ \ \ \  + (\sigma_{3}^{l})^{1-p}(\sigma_{3}^{u})^{p}y(t)dB_{3}(t) +  \int_{\mathbb{Y}}c_{3}(u)y(t^{-}){N}(dt,du),
\end{array}  
\right.
\end{equation}
for $ p \in [0, 1] $.

Throughout this paper, let 
\[
b_{i} = (r_{i}^{l})^{1-p}(r_{i}^{u})^{p} - \dfrac{{((\sigma_{i}^{l})^{1-p}(\sigma_{i}^{u})^{p})}^{2}}{2} + \int_{\mathbb{Y}}  \ln (1+ c_{i}(u) )\lambda (du), \  i=1,2.
\]

\[
b_{3} = -(r_{3}^{l})^{1-p}(r_{3}^{u})^{p} - \dfrac{{((\sigma_{3}^{l})^{1-p}(\sigma_{3}^{u})^{p})}^{2}}{2} + \int_{\mathbb{Y}}  \ln (1+ c_{3}(u) )\lambda (du).
\]

\[
\langle f(t) \rangle = t^{-1} \int_{0}^{t} f(s) ds, \ \  \langle f(t) \rangle^{*} = \limsup_{t\rightarrow \infty} t^{-1} \int_{0}^{t} f(s) ds, \ \  \langle f(t) \rangle_{*} = \liminf_{t\rightarrow \infty} t^{-1} \int_{0}^{t} f(s) ds.
\]

\section{ Existence and uniqueness of positive solution of system (\ref{2.3}) }
\label{sec:3}

For convenience in the following investigation, we  require that
\begin{description}
	\item[(H1)] $ \int_{\mathbb{Y}}  ( |c_{i}(u)|\vee |c_{i}(u)|^{2}   )\lambda (du) \leq \infty$.
	\item[(H2)] $\int_{\mathbb{Y}}  ( | \ln (1+ c_{i}(u) )|\vee | \ln (1+ c_{i}(u) )|^{2}   )\lambda (du) \leq \infty \ \ i=1,2,3$ .
\end{description}

The following theorem will prove that system (\ref{2.3}) admits a unique global positive solution.

\begin{Theorem}\label{T31}
	Let Assumptions (H1) and (H2) hold. Then for any given initial value $ (x_{1}(0), x_{2}(0), y(0)) \in \mathbb{R}^{3}_{+} $, system (\ref{2.3}) has a unique solution $ (x_{1}(t), x_{2}(t), y(t)) \in \mathbb{R}^{3}_{+} $ for all $ t \geq 0 $ almost surely (a.s.).
\end{Theorem}

\noindent \textbf{Proof.} 
Since $ t \geq 0 $, by system (\ref{2.3}), we can construct the following system
\begin{equation*}
\left\{
\begin{array}{lcl}
du_{1}(t) =  \Big(b_{1} - (a_{11}^{l})^{1-p}(a_{11}^{u})^{p} e^{u_{1}(t)} - (a_{12}^{l})^{1-p}(a_{12}^{u})^{p}e^{u_{2}(t)} - \dfrac{(a_{13}^{l})^{1-p}(a_{13}^{u})^{p}e^{v(t)}}{1+ e^{u_{1}(t)}}\Big)dt \\
 \ \ \ \ \ \ \ \ \ \ \ \ + (\sigma_{1}^{l})^{1-p}(\sigma_{1}^{u})^{p}dB_{1}(t) +  \int_{\mathbb{Y}}\ln(1+c_{1}(u)) \widetilde{N}(dt, du),\\
du_{2}(t) =  \Big(b_{2} - (a_{21}^{l})^{1-p}(a_{21}^{u})^{p} e^{u_{1}(t)} - (a_{22}^{l})^{1-p}(a_{22}^{u})^{p}e^{u_{2}(t)} - \dfrac{(a_{23}^{l})^{1-p}(a_{23}^{u})^{p}e^{v(t)}}{1+ e^{u_{2}(t)}}\Big)dt\\
 \ \ \ \ \ \ \ \ \ \ \ \  + (\sigma_{2}^{l})^{1-p}(\sigma_{2}^{u})^{p}dB_{2}(t) +  \int_{\mathbb{Y}}\ln(1+c_{2}(u)) \widetilde{N}(dt, du),\\
dv(t) =  \Big(b_{3} - (a_{33}^{l})^{1-p}(a_{33}^{u})^{p} e^{v(t)} + \dfrac{(a_{31}^{l})^{1-p}(a_{31}^{u})^{p}e^{u_{1}(t)}}{1+ e^{u_{1}(t)}} + \dfrac{(a_{32}^{l})^{1-p}(a_{32}^{u})^{p}e^{u_{2}(t)}}{1+ e^{u_{2}(t)}}\Big)dt\\
 \ \ \ \ \ \ \ \ \ \ \ \  + (\sigma_{3}^{l})^{1-p}(\sigma_{3}^{u})^{p}dB_{3}(t) +  \int_{\mathbb{Y}}\ln(1+c_{3}(u)) \widetilde{N}(dt, du),
\end{array}  
\right.
\end{equation*}
Because the coefficients of this system are local Lipschitz continuous (Mao \cite{25}), for any given initial value $ (u_{1}(0), u_{2}(0), v(0))= ( \ln x_{1}(0), \ln x_{2}(0), \ln y(0) ) \in \mathbb{R}^{3}_{+} $, there is a unique local solution $ (u_{1}(t), u_{2}(t), v(t)) $ on $ t\in [0,\tau) $, where $ \tau $ is the explosion time (see Mao \cite{25}). Hence, system (\ref{2.3}) admits unique positive local solution $ (x_{1}(t), x_{2}(t), y(t)) = ( e^{u_{1}(t)}, e^{u_{2}(t)}, e^{v(t)} ) $. The proof of global existence of this local solution to system (\ref{2.3}) is rather standard. Therefore, we omit the proof here.\qed\\

By Theorem~\ref{T31}, system (\ref{2.3}) admits a unique global positive solution. Next, we will show that the solution of system (\ref{2.3}) is stochastically bounded.

\begin{Theorem}
	Let Assumptions (H1) and (H2) hold. Then for any given initial value $ (x_{1}(0), x_{2}(0), y(0)) \in \mathbb{R}^{3}_{+} $ and $ k>0 $, the solution of system (\ref{2.3}) satisfies
	\[
	\limsup_{t\rightarrow \infty} \mathbb{E} \big( x_{1}^{k}(t) + x_{2}^{k}(t) + y^{k}(t) \big) \leq K,
	\]
where $ K $ is a generic positive constant.
\end{Theorem}

\noindent \textbf{Proof.} 
The It\v{o}’s formula (Situ \cite{26}) shows that 
\[
\mathbb{E} \big( e^{t} (x_{1}^{k}(t) + x_{2}^{k}(t) + y^{k}(t)) \big) =x_{1}^{k}(0) + x_{2}^{k}(0) + y^{k}(0) + \mathbb{E}\int_{0}^{t} e^{s} F(s)ds,
\]
where
\begin{equation*}
\begin{split}
F =& -{(a_{11}^{l})^{1-p}(a_{11}^{u})^{p}}k x_{1}^{k+1} + \Big( 1+ k{(r_{1}^{l})^{1-p}(r_{1}^{u})^{p}} + \dfrac{k(k-1)}{2}({(\sigma_{1}^{l})^{1-p}(\sigma_{1}^{u})^{p}})^{2} + \int_{\mathbb{Y}}\big( (1+c_{1}(u))^{k} -1\big)\lambda (du)    \Big) x_{1}^{k}\\
& - {(a_{12}^{l})^{1-p}(a_{12}^{u})^{p}}kx_{1}^{k}x_{2} - \dfrac{{(a_{13}^{l})^{1-p}(a_{13}^{u})^{p}}kx_{1}^{k}y}{1+x_{1}}\\
&-{(a_{22}^{l})^{1-p}(a_{22}^{u})^{p}}k x_{2}^{k+1} + \Big( 1+ k{(r_{2}^{l})^{1-p}(r_{2}^{u})^{p}} + \dfrac{k(k-1)}{2}({(\sigma_{2}^{l})^{1-p}(\sigma_{2}^{u})^{p}})^{2} + \int_{\mathbb{Y}}\big( (1+c_{2}(u))^{k} -1\big)\lambda (du)    \Big) x_{2}^{k}\\
& - {(a_{21}^{l})^{1-p}(a_{21}^{u})^{p}}kx_{1}x_{2}^{k} - \dfrac{{(a_{23}^{l})^{1-p}(a_{23}^{u})^{p}}kx_{2}^{k}y}{1+x_{2}}\\
&-{(a_{33}^{l})^{1-p}(a_{33}^{u})^{p}}k y^{k+1} + \Big( 1- k{(r_{3}^{l})^{1-p}(r_{3}^{u})^{p}} + \dfrac{k(k-1)}{2}({(\sigma_{3}^{l})^{1-p}(\sigma_{3}^{u})^{p}})^{2} + \int_{\mathbb{Y}}\big( (1+c_{3}(u))^{k} -1\big)\lambda (du)    \Big) y^{k} \\
&+ \dfrac{{(a_{31}^{l})^{1-p}(a_{31}^{u})^{p}}kx_{1}y^{k}}{1+x_{1}} + \dfrac{{(a_{32}^{l})^{1-p}(a_{32}^{u})^{p}}kx_{2}y^{k}}{1+x_{2}}\\
\leq &  -{(a_{11}^{l})^{1-p}(a_{11}^{u})^{p}}k x_{1}^{k+1} + \Big( 1+ k{(r_{1}^{l})^{1-p}(r_{1}^{u})^{p}} + \dfrac{k(k-1)}{2}({(\sigma_{1}^{l})^{1-p}(\sigma_{1}^{u})^{p}})^{2} + \int_{\mathbb{Y}}\big( (1+c_{1}(u))^{k} -1\big)\lambda (du)    \Big) x_{1}^{k}\\
&-{(a_{22}^{l})^{1-p}(a_{22}^{u})^{p}}k x_{2}^{k+1} + \Big( 1+ k{(r_{2}^{l})^{1-p}(r_{2}^{u})^{p}} + \dfrac{k(k-1)}{2}({(\sigma_{2}^{l})^{1-p}(\sigma_{2}^{u})^{p}})^{2} + \int_{\mathbb{Y}}\big( (1+c_{2}(u))^{k} -1\big)\lambda (du)    \Big) x_{2}^{k}\\
&-{(a_{33}^{l})^{1-p}(a_{33}^{u})^{p}}k y^{k+1} + \Big( 1- k{(r_{3}^{l})^{1-p}(r_{3}^{u})^{p}} + {(a_{31}^{l})^{1-p}(a_{31}^{u})^{p}}k + {(a_{32}^{l})^{1-p}(a_{32}^{u})^{p}}k \\
&+ \dfrac{k(k-1)}{2}({(\sigma_{3}^{l})^{1-p}(\sigma_{3}^{u})^{p}})^{2} + \int_{\mathbb{Y}}\big( (1+c_{3}(u))^{k} -1\big)\lambda (du)    \Big) y^{k}\\
\leq & K.
\end{split}
\end{equation*}
Thus,
\[
\mathbb{E} \big( e^{t} (x_{1}^{k}(t) + x_{2}^{k}(t) + y^{k}(t)) \big) \leq x_{1}^{k}(0) + x_{2}^{k}(0) + y^{k}(0) + K e^{t},
\]
which implies
	\[
\limsup_{t\rightarrow \infty} \mathbb{E} \big( x_{1}^{k}(t) + x_{2}^{k}(t) + y^{k}(t) \big) \leq K.
\]
This completes the proof.  \qed

\section{Extinction}
\label{sec:4}

When studying mathematical ecology, two of the most interesting issues are persistence and extinction. In this section, we discuss the extinction of populations in system (\ref{2.3}) and leave its persistence to the next section.\\

\begin{Definition}
	\item[(1)] 	$ x(t) $ is said to be extinct if $ \lim_{t\rightarrow \infty}x(t) =0 $ a.s.
	\item[(2)] 	$ x(t) $ is said to be strongly persistent in the mean if $ \liminf_{t\rightarrow \infty} t^{-1} \int_{0}^{t} x(s) ds > 0 $ a.s.
\end{Definition}

Before we state the main results of this section, several lemmas (Zhang et al. \cite{27}) should be recalled without the proofs and relevant explanations.\\

\begin{Lemma}\label{L41}
For any given initial value $ (x_{1}(0), x_{2}(0), y(0)) \in \mathbb{R}^{3}_{+} $, the solution $ (x_{1}(t), x_{2}(t), y(t))$ of system (\ref{2.3}) satisfies 
\[
\limsup_{t\rightarrow \infty} \dfrac{\ln x_{i}(t)}{t} \leq 0, \ \ \ \  \limsup_{t\rightarrow \infty} \dfrac{\ln y(t)}{t} \leq 0 \ \ a.s. \ \ \ \ \ \ i=1,2.
\]
\end{Lemma}

\begin{Lemma}\label{L42}
Assume
\[
dx(t) = x(t) \big((r^{l})^{1-p}(r^{u})^{p} - (a^{l})^{1-p}(a^{u})^{p}x(t)\big)dt + (\sigma^{l})^{1-p}(\sigma^{u})^{p}x(t)dB(t) +  \int_{\mathbb{Y}}c(u)x(t^{-}){N}(dt,du).
\]	
If $ (a^{l})^{1-p}(a^{u})^{p}>0 $	and $ (r^{l})^{1-p}(r^{u})^{p} - \dfrac{({ (\sigma^{l})^{1-p}(\sigma^{u})^{p} })^{2}}{2} + \int_{\mathbb{Y}}\ln (1+c(u)) \lambda (du) \geq 0$, we have 
\[
\lim_{t\rightarrow \infty} \dfrac{\ln x(t)}{t} =0 \ a.s.
\]
and
\[
\lim_{t\rightarrow \infty} \dfrac{1}{t} \int_{0}^{t} x(s) ds = \dfrac{(r^{l})^{1-p}(r^{u})^{p} - \dfrac{({ (\sigma^{l})^{1-p}(\sigma^{u})^{p} })^{2}}{2} + \int_{\mathbb{Y}}\ln (1+c(u)) \lambda (du)}{(a^{l})^{1-p}(a^{u})^{p}}  \  a.s.
\]
\end{Lemma}

We now establish sufficient conditions for extinction of populations in system (\ref{2.3}).

\begin{Theorem}\label{T41}
	If $ b_{1}>0 $, $ b_{2}>0 $ and $ b_{3}<0 $, all the populations in system (\ref{2.3}) go to extinction.
\end{Theorem}

\noindent \textbf{Proof.}
The It\v{o}’s formula of system (\ref{2.3}) yields 
\begin{equation}\label{41}
\begin{split}
d\ln x_{1}(t) =&  \Big(b_{1} - (a_{11}^{l})^{1-p}(a_{11}^{u})^{p}x_{1}(t) - (a_{12}^{l})^{1-p}(a_{12}^{u})^{p}x_{2}(t) - \dfrac{(a_{13}^{l})^{1-p}(a_{13}^{u})^{p}y(t)}{1+x_{1}(t)}\Big)dt \\
+& (\sigma_{1}^{l})^{1-p}(\sigma_{1}^{u})^{p}dB_{1}(t) +  \int_{\mathbb{Y}} \ln (1+c_{1}(u))\widetilde{N}(dt,du)
\end{split}
\end{equation}
(\ref{41}) from $ 0 $ to $ t $ and then dividing by $ t $ on both sides, we obtain
\begin{equation}\label{42}
\begin{split}
\dfrac{\ln(x_{1}(t)/x_{1}(0))}{t} =& b_{1} - (a_{11}^{l})^{1-p}(a_{11}^{u})^{p}  \langle x_{1}(t) \rangle - (a_{12}^{l})^{1-p}(a_{12}^{u})^{p} \langle x_{2}(t) \rangle - (a_{13}^{l})^{1-p}(a_{13}^{u})^{p}\langle \dfrac{y(t)}{1+x_{1}(t)}\rangle\\
& + \dfrac{M_{1}(t)}{t} + \dfrac{\widetilde{M}_{1}(t)}{t}.
\end{split}
\end{equation}
Similarly, we have
\begin{equation}\label{43}
\begin{split}
\dfrac{\ln(x_{2}(t)/x_{2}(0))}{t} =& b_{2} - (a_{21}^{l})^{1-p}(a_{21}^{u})^{p}  \langle x_{1}(t) \rangle - (a_{22}^{l})^{1-p}(a_{22}^{u})^{p} \langle x_{2}(t) \rangle - (a_{23}^{l})^{1-p}(a_{23}^{u})^{p}\langle \dfrac{y(t)}{1+x_{2}(t)}\rangle\\
 &+ \dfrac{M_{2}(t)}{t} + \dfrac{\widetilde{M}_{2}(t)}{t},
\end{split}
\end{equation}
and
\begin{equation}\label{44}
\begin{split}
\dfrac{\ln(y(t)/y(0))}{t} =& b_{3} - (a_{33}^{l})^{1-p}(a_{33}^{u})^{p}  \langle y(t) \rangle + (a_{31}^{l})^{1-p}(a_{31}^{u})^{p}\langle \dfrac{x_{1}(t)}{1+x_{1}(t)}\rangle + (a_{32}^{l})^{1-p}(a_{32}^{u})^{p}\langle \dfrac{x_{2}(t)}{1+x_{2}(t)}\rangle\\ &+ \dfrac{M_{3}(t)}{t} + \dfrac{\widetilde{M}_{3}(t)}{t}.
\end{split}
\end{equation}
where $ M_{i}(t): = \int_{0}^{t}{(\sigma_{i}^{l})^{1-p} (\sigma_{i}^{u})^{p}}dB_{i}(s) $ and $ \widetilde{M}_{i}(t):= \int_{0}^{t}\int_{\mathbb{Y}}[\ln (1+c_{i}(u) )] \widetilde{N}(ds,du)   $, $ i=1,2,3 $ are all martingale terms. Thus, by strong law of large numbers, we have
\begin{equation}\label{45}
\lim_{t\rightarrow \infty} \dfrac{M_{i}(t)}{t} = 0 \ \ a.s. \textrm{ and } \lim_{t\rightarrow \infty} \dfrac{\widetilde{M}_{i}(t)}{t} = 0 \ \ a.s.
\end{equation}
Thus, (\ref{42}), (\ref{43}) and (\ref{45}) yields 
\[
\limsup_{t\rightarrow \infty} \dfrac{\ln x_{i}(t)}{t} \leq b_{i} \ \ a.s. \ \ \ \ \ i=1,2.
\]
Which means $ \lim_{t\rightarrow \infty} x_{i}(t) =0 $ a.s. $ i=1,2 $ when $ b_{1}<0 $ and $ b_{2}<0 $. This together with (\ref{44}) and (\ref{45}) implies
\[
\limsup_{t\rightarrow \infty} \dfrac{\ln y(t)}{t} \leq b_{3} \ \ a.s. 
\]
Therefore, $ b_{3} < 0 $ implies $ \lim_{t\rightarrow \infty} y(t) =0 $ a.s. This completes the proof.  \qed

\begin{Theorem}
	If $ b_{1}>0 $, $ b_{2}<0 $ and $ b_{3} + (a_{31}^{l})^{1-p}(a_{31}^{u})^{p} <0$, $ x_{2}(t) $ and $ y(t) $ are extinct and
	\[
	\limsup_{t\rightarrow \infty} \dfrac{1}{t} \int_{0}^{t} x_{1}(s)ds = \dfrac{b_{1}}{(a_{11}^{l})^{1-p}(a_{11}^{u})^{p}} \ \ a.s.
	\]
\end{Theorem}

\noindent \textbf{Proof.}
$ b_{2}<0 $ toghter with (\ref{43}) and (\ref{45}) yields
\[
\limsup_{t\rightarrow \infty}\dfrac{\ln x_{2}(t)}{t} \leq b_{2} <0 \ \ a.s.,
\]
which means $ \lim_{t\rightarrow \infty} x_{2}(t) = 0 $. Combining this with (\ref{44}) and (\ref{45}), we know that
\[
\limsup_{t\rightarrow \infty}\dfrac{\ln y(t)}{t} \leq b_{3} + (a_{31}^{l})^{1-p}(a_{31}^{u})^{p} <0 \ \ a.s.,
\]
which also implies  $ \lim_{t\rightarrow \infty} y(t) = 0 $.
It is easy to check that, for any $ 0< \epsilon < \dfrac{b_{1}}{2} $, there exist a positive constant $ t_{0} $ and a
set $ \Omega_{\epsilon} $ such that $ \mathbb{P}(\Omega_{\epsilon}) \geq 1-\epsilon $, and for $ t\geq t_{0} $ we get $ (a_{12}^{l})^{1-p}(a_{12}^{u})^{p}x_{2} < \epsilon $, $ (a_{13}^{l})^{1-p}(a_{13}^{u})^{p}y < \epsilon $. Thus, for any $ \omega\in \Omega_{\epsilon} $,
\[
dx_{1}(t) \leq x_{1}(t) [(r_{1}^{l})^{1-p}(r_{1}^{u})^{p} - (a_{11}^{l})^{1-p}(a_{11}^{u})^{p}x_{1}(t)]dt + (\sigma_{1}^{l})^{1-p}(\sigma_{1}^{u})^{p}x_{1}(t)dB_{1}(t) +  \int_{\mathbb{Y}}c_{1}(u)x_{1}(t^{-}){N}(dt,du),
\]
and 
\[
dx_{1}(t) \geq x_{1}(t) [(r_{1}^{l})^{1-p}(r_{1}^{u})^{p} -2\epsilon - (a_{11}^{l})^{1-p}(a_{11}^{u})^{p}x_{1}(t)]dt + (\sigma_{1}^{l})^{1-p}(\sigma_{1}^{u})^{p}x_{1}(t)dB_{1}(t) +  \int_{\mathbb{Y}}c_{1}(u)x_{1}(t^{-}){N}(dt,du).
\]
By Lemma~\ref{L42} and stochastic comparison theorem, for $ b_{1} > 0 $, 
\[
\dfrac{b_{1}- 2\epsilon}{(a_{11}^{l})^{1-p}(a_{11}^{u})^{p}} \leq \liminf_{t\rightarrow \infty} \dfrac{1}{t} \int_{0}^{t} x_{1}(s)ds \leq \limsup_{t\rightarrow \infty} \dfrac{1}{t} \int_{0}^{t} x_{1}(s)ds \leq \dfrac{b_{1}}{(a_{11}^{l})^{1-p}(a_{11}^{u})^{p}} \ \ a.s.
\]
Thus, we have 
\[
\lim_{t\rightarrow \infty} \dfrac{1}{t} \int_{0}^{t} x_{1}(s)ds = \dfrac{b_{1}}{(a_{11}^{l})^{1-p}(a_{11}^{u})^{p}} \ \ a.s.,
\]
when $ \epsilon \rightarrow 0 $.  This completes the proof.  \qed

Next, we establish sufficient conditions for persistence in the mean of system (\ref{2.3}). Before that, we need to consider several stochastic differential equations with jumps.

\[d\phi_{1}(t) = \phi_{1}(t) [(r_{1}^{l})^{1-p}(r_{1}^{u})^{p} - (a_{11}^{l})^{1-p}(a_{11}^{u})^{p}\phi_{1}(t)]dt + (\sigma_{1}^{l})^{1-p}(\sigma_{1}^{u})^{p}\phi_{1}(t)dB_{1}(t) +  \int_{\mathbb{Y}}c_{1}(u)\phi_{1}(t^{-}){N}(dt,du),\]
\[d\phi_{2}(t) = \phi_{2}(t) [(r_{2}^{l})^{1-p}(r_{2}^{u})^{p} -  (a_{22}^{l})^{1-p}(a_{22}^{u})^{p}\phi_{2}(t)]dt + (\sigma_{2}^{l})^{1-p}(\sigma_{2}^{u})^{p}\phi_{2}(t)dB_{2}(t) +  \int_{\mathbb{Y}}c_{2}(u)\phi_{2}(t^{-}){N}(dt,du),\]
\[d\phi_{3}(t) = \phi_{3}(t) [-(r_{3}^{l})^{1-p}(r_{3}^{u})^{p} -  (a_{33}^{l})^{1-p}(a_{33}^{u})^{p}\phi_{3}(t)]dt + (\sigma_{3}^{l})^{1-p}(\sigma_{3}^{u})^{p}\phi_{3}(t)dB_{3}(t) +  \int_{\mathbb{Y}}c_{3}(u)\phi_{3}(t^{-}){N}(dt,du),\]
\begin{equation*}
\begin{split}
d\phi_{4}(t) =& \phi_{4}(t) [-(r_{3}^{l})^{1-p}(r_{3}^{u})^{p} + (a_{31}^{l})^{1-p}(a_{31}^{u})^{p} + (a_{32}^{l})^{1-p}(a_{32}^{u})^{p} -  (a_{33}^{l})^{1-p}(a_{33}^{u})^{p}\phi_{4}(t)]dt \\ &+ (\sigma_{3}^{l})^{1-p}(\sigma_{3}^{u})^{p}\phi_{4}(t)dB_{3}(t) +  \int_{\mathbb{Y}}c_{3}(u)\phi_{4}(t^{-}){N}(dt,du).
\end{split}
\end{equation*}
 Thus, we know from stochastic comparison theorem that 
 \[
 x_{1}(t) \leq \phi_{1}(t), \ \ \  x_{2}(t) \leq \phi_{2}(t), \ \ \   \phi_{3}(t) \leq y(t) \leq \phi_{4}(t).
 \]
 
 \begin{Theorem}
 	If $ b_{3}>0 $,  population $ y(t) $ of system (\ref{2.3}) satisfies
 	\[
 	\dfrac{b_{3}}{(a_{33}^{l})^{1-p}(a_{33}^{u})^{p}} \leq \liminf_{t\rightarrow \infty} \dfrac{1}{t} \int_{0}^{t} y(s)ds \leq \limsup_{t\rightarrow \infty} \dfrac{1}{t} \int_{0}^{t} y(s)ds \leq \dfrac{b_{3} + (a_{31}^{l})^{1-p}(a_{31}^{u})^{p} + (a_{32}^{l})^{1-p}(a_{32}^{u})^{p}}{(a_{33}^{l})^{1-p}(a_{33}^{u})^{p}} \ \ a.s.
 	\]
 	Moveover, if \\$ b_{1}> \max\{ 0, ((a_{12}^{l})^{1-p}(a_{12}^{u})^{p})\dfrac{b_{2}}{(a_{22}^{l})^{1-p}(a_{22}^{u})^{p}} + ((a_{13}^{l})^{1-p}(a_{13}^{u})^{p})\dfrac{b_{3} + (a_{31}^{l})^{1-p}(a_{31}^{u})^{p} + (a_{32}^{l})^{1-p}(a_{32}^{u})^{p}}{(a_{33}^{l})^{1-p}(a_{33}^{u})^{p}} \} $\\ and\\ $ b_{2}> \max\{ 0, ((a_{21}^{l})^{1-p}(a_{21}^{u})^{p})\dfrac{b_{1}}{(a_{11}^{l})^{1-p}(a_{11}^{u})^{p}} + ((a_{23}^{l})^{1-p}(a_{23}^{u})^{p})\dfrac{b_{3} + (a_{31}^{l})^{1-p}(a_{31}^{u})^{p} + (a_{32}^{l})^{1-p}(a_{32}^{u})^{p}}{(a_{33}^{l})^{1-p}(a_{33}^{u})^{p}} \} $,\\ we have 
 	\begin{equation*}
 	\begin{split}
 	\dfrac{1}{(a_{11}^{l})^{1-p}(a_{11}^{u})^{p}}\Big( b_{1} -&  ((a_{12}^{l})^{1-p}(a_{12}^{u})^{p})\dfrac{b_{2}}{(a_{22}^{l})^{1-p}(a_{22}^{u})^{p}} - ((a_{13}^{l})^{1-p}(a_{13}^{u})^{p})\dfrac{b_{3} + (a_{31}^{l})^{1-p}(a_{31}^{u})^{p} + (a_{32}^{l})^{1-p}(a_{32}^{u})^{p}}{(a_{33}^{l})^{1-p}(a_{33}^{u})^{p}} \Big)\\
 	& \leq \liminf_{t\rightarrow \infty} \dfrac{1}{t} \int_{0}^{t} x_{1}(s)ds \leq \limsup_{t\rightarrow \infty} \dfrac{1}{t} \int_{0}^{t} x_{1}(s)ds \leq \dfrac{b_{1}}{(a_{11}^{l})^{1-p}(a_{11}^{u})^{p}}\ \ a.s.,
 \end{split}
 \end{equation*}
 	and 
 		\begin{equation*}
 	\begin{split}
 	\dfrac{1}{(a_{22}^{l})^{1-p}(a_{22}^{u})^{p}}\Big( b_{2} -&  ((a_{21}^{l})^{1-p}(a_{21}^{u})^{p})\dfrac{b_{1}}{(a_{11}^{l})^{1-p}(a_{11}^{u})^{p}} - ((a_{23}^{l})^{1-p}(a_{23}^{u})^{p})\dfrac{b_{3} + (a_{31}^{l})^{1-p}(a_{31}^{u})^{p} + (a_{32}^{l})^{1-p}(a_{32}^{u})^{p}}{(a_{33}^{l})^{1-p}(a_{33}^{u})^{p}} \Big)\\
 	& \leq \liminf_{t\rightarrow \infty} \dfrac{1}{t} \int_{0}^{t} x_{2}(s)ds \leq \limsup_{t\rightarrow \infty} \dfrac{1}{t} \int_{0}^{t} x_{2}(s)ds \leq \dfrac{b_{2}}{(a_{22}^{l})^{1-p}(a_{22}^{u})^{p}}\ \ a.s.,
 	\end{split}
 	\end{equation*}
 	which means all the  populations in system (\ref{2.3}) are strongly persistent in the mean.
 \end{Theorem}

\noindent \textbf{Proof.}
By
\[
\phi_{3}(t) \leq y(t) \leq \phi_{4}(t),
\]
we know from Lemma~\ref{L42} that if $ b_{3}>0 $,
\[
\dfrac{b_{3}}{(a_{33}^{l})^{1-p}(a_{33}^{u})^{p}} \leq \liminf_{t\rightarrow \infty} \dfrac{1}{t} \int_{0}^{t} y(s)ds \leq \limsup_{t\rightarrow \infty} \dfrac{1}{t} \int_{0}^{t} y(s)ds \leq \dfrac{b_{3} + (a_{31}^{l})^{1-p}(a_{31}^{u})^{p} + (a_{32}^{l})^{1-p}(a_{32}^{u})^{p}}{(a_{33}^{l})^{1-p}(a_{33}^{u})^{p}} \ \ a.s.
\]
By Lemma~\ref{L42} we still know that 
\[
\limsup_{t\rightarrow \infty} \dfrac{1}{t} \int_{0}^{t} x_{1}(s)ds  \leq \lim_{t\rightarrow \infty} \dfrac{1}{t} \int_{0}^{t} \phi_{1}(s)ds = \dfrac{b_{1}}{(a_{11}^{l})^{1-p}(a_{11}^{u})^{p}}\ \ a.s.,
\]
and 
\[
\limsup_{t\rightarrow \infty} \dfrac{1}{t} \int_{0}^{t} x_{2}(s)ds  \leq \lim_{t\rightarrow \infty} \dfrac{1}{t} \int_{0}^{t} \phi_{2}(s)ds = \dfrac{b_{2}}{(a_{22}^{l})^{1-p}(a_{22}^{u})^{p}}\ \ a.s.
\]
These together with Lemma~\ref{L41} and (\ref{42})  (\ref{45}) implies
\begin{equation*}
\begin{split}
(a_{11}^{l})^{1-p}(a_{11}^{u})^{p} \liminf_{t\rightarrow \infty} \dfrac{1}{t} \int_{0}^{t} x_{1}(s)ds  \geq & \liminf_{t\rightarrow \infty} \Big\{ -\dfrac{\ln(x_{1}(t)/x_{1}(0))}{t} + b_{1}  - (a_{12}^{l})^{1-p}(a_{12}^{u})^{p} \dfrac{1}{t}\int_{0}^{t} x_{2}(s)ds\\
&-(a_{13}^{l})^{1-p}(a_{13}^{u})^{p} \dfrac{1}{t} \int_{0}^{t} \dfrac{y(s)}{1+x(s)}ds +  \dfrac{M_{1}(t)}{t} + \dfrac{\widetilde{M}_{1}(t)}{t}  \Big\}\\
\geq & b_{1} - \limsup_{t\rightarrow \infty} \dfrac{\ln x_{1}(t)}{t} - (a_{12}^{l})^{1-p}(a_{12}^{u})^{p} \limsup_{t\rightarrow \infty} \dfrac{1}{t} \int_{0}^{t} \phi_{2}(s)ds \\
&- (a_{13}^{l})^{1-p}(a_{13}^{u})^{p} \limsup_{t\rightarrow \infty} \dfrac{1}{t} \int_{0}^{t} y(s)ds\\
\geq & b_{1} - (a_{12}^{l})^{1-p}(a_{12}^{u})^{p}\dfrac{b_{2}}{(a_{22}^{l})^{1-p}(a_{22}^{u})^{p}}\\
& -  (a_{13}^{l})^{1-p}(a_{13}^{u})^{p} \dfrac{b_{3} + (a_{31}^{l})^{1-p}(a_{31}^{u})^{p} + (a_{32}^{l})^{1-p}(a_{32}^{u})^{p}}{(a_{33}^{l})^{1-p}(a_{33}^{u})^{p}} \ \ a.s.
\end{split}
\end{equation*}
Similarly, we also have
\begin{equation*}
\begin{split}
(a_{22}^{l})^{1-p}(a_{22}^{u})^{p} \liminf_{t\rightarrow \infty} \dfrac{1}{t} \int_{0}^{t} x_{2}(s)ds \geq & b_{2} - (a_{21}^{l})^{1-p}(a_{21}^{u})^{p}\dfrac{b_{1}}{(a_{11}^{l})^{1-p}(a_{11}^{u})^{p}}\\
& -  (a_{23}^{l})^{1-p}(a_{23}^{u})^{p} \dfrac{b_{3} + (a_{31}^{l})^{1-p}(a_{31}^{u})^{p} + (a_{32}^{l})^{1-p}(a_{32}^{u})^{p}}{(a_{33}^{l})^{1-p}(a_{33}^{u})^{p}} \ \ a.s.
\end{split}
\end{equation*}
 This completes the proof.  \qed

\end{document}